\documentclass[10pt,letterpaper,conference]{IEEEconf}  
 \usepackage{etoolbox}
\makeatletter
\patchcmd{\@begintheorem}{\textit}{\textbf}{}{}
\makeatother
 \newtheorem{definition}{\bf Definition}
 
  \newtheorem{Corollary}{\bf Corollary}
  \newtheorem{thm}{\bf Theorem}

 \newtheorem{lemma}{\bf Lemma}
  \newtheorem{prop}{\bf Proposition}
  \usepackage{dirtytalk}
\usepackage[]{amsmath}
\usepackage{amsfonts}
\usepackage{amssymb}
\usepackage{bbm}
\usepackage{breqn}
\usepackage{subcaption}
\usepackage{cite}
\usepackage{cases}
\usepackage{tabu}
\usepackage{url}
\usepackage{multicol}
\usepackage{color}
\usepackage[bottom]{footmisc}
\usepackage{environ}
\usepackage{tikz} 
\usepackage{enumerate}
\usetikzlibrary{arrows,decorations.markings}
\usetikzlibrary{decorations.pathreplacing}

\usetikzlibrary{automata,arrows,positioning,calc}
\usepackage[]{algpseudocode}
\algtext*{EndWhile}
\algtext*{EndIf}
\algtext*{EndFor}
\usepackage{algorithm}
\title{\LARGE \bf Incentive Design for Temporal Logic Objectives}
\IEEEoverridecommandlockouts
\author{Yagiz Savas, Vijay Gupta, Melkior Ornik, Lillian J. Ratliff, Ufuk Topcu \thanks{ Y. Savas and U. Topcu are with the Department of Aerospace Engineering, University of Texas at Austin, TX, USA. E-mail: \{yagiz.savas, utopcu\}@utexas.edu   } \thanks{V. Gupta is with the Department of Electrical Engineering, University of Notre Dame, IN, USA. E-mail: vgupta@nd.edu}  \thanks{M. Ornik is with the Department of Aerospace Engineering and the Coordinated Science Laboratory, University of Illinois at Urbana-Champaign, IL, USA. E-mail: mornik@illinois.edu}  \thanks{L. J. Ratliff is with the Department of Electrical Engineering, University of Washington, WA, USA. E-mail: ratliffl@uw.edu}  }
\date{}
\begin{document}
\maketitle
{\abstract  We study the problem of designing an optimal sequence of incentives that a principal should offer to an agent so that the agent's optimal behavior under the incentives realizes the principal's objective expressed as a temporal logic formula. We consider an agent with a finite decision horizon and model its decision-making process as a Markov decision process (MDP). Under certain assumptions, we present a polynomial-time algorithm to synthesize an incentive sequence that minimizes the cost to the principal. We show that if the underlying MDP has only deterministic transitions, the principal can hide its objective from the agent and still realize the desired behavior through incentives. On the other hand, an MDP with stochastic transitions may require the principal to share its objective with the agent. Finally, we demonstrate the proposed method in motion planning examples where a principal changes the optimal trajectory of an agent by providing incentives.  

}

\section{Introduction}
Consider a scenario where a principal provides incentives to an agent so that the optimal behavior of the agent under the provided incentives satisfies the principal's objective. If the principal had enough resources to provide arbitrarily large incentives, it would be straightforward to obtain the desired agent behaviour. However, since the resources are limited in practice, it is important to establish the minimum amount of incentives that leads to the desired behavior. In this paper, we are interested in designing a sequence of incentives that minimizes the cost to the principal while guaranteeing the realization of its objective by the agent with maximum probability.

We model the sequential decision-making process of the agent as a Markov decision process (MDP) \cite{Puterman}, and assume that the agent's objective is to maximize its expected total reward at the end of a finite planning horizon. Although each planning horizon is finite, the agent plans its future decisions infinitely many times. Examples of such an agent can be a person who plans her schedule on a weekly basis or an autonomous system with a limited computational power which plans its route by considering only a small subset of all possible environment states. 

The principal's objective is described by a syntactically co-safe linear temporal logic (LTL) formula. LTL specifications are widely used to describe complex tasks for autonomous robots \cite{Kloetzer}, design security protocols \cite{Armando} and check the reliability of software \cite{Tan}. For example, in a navigation scenario, syntactically co-safe LTL formulae allow one to specify tasks such as liveness (eventually visit the region A) or priority (first visit the region A and then B).

We assume that the principal is aware of the agent's reward function and the length of its planning horizon. In many real-world applications, the decision horizon and the reward structure of an agent can be known or at least inferred through observations. For example, a manufacturing company is generally interested in maximizing its profit at the end of a fiscal year, and an autonomous car aims to reach its destination within certain time interval. 

From a practical point of view, an interesting question is whether an adversarial principal can convince an agent to satisfy its objective through incentives. In such a scenario, if the agent knows the principal's objective explicitly, it will reject the provided incentives because the resulting behavior under the incentives will serve to the benefit of the enemy. However, if the principal can design an incentive sequence without sharing its objective with the agent, then the incentives may lead to the desired agent behavior. Therefore, it is important to establish the conditions under which the principal can actually hide its objective from the agent.

The contributions of this paper can be summarized as follows. First, we present an algorithm, based on a series of linear optimization problems, to synthesize a sequence of incentives that minimizes the cost to the principal while ensuring that the optimal agent behavior under the provided incentives satisfies a syntactically co-safe LTL formula with maximum probability. Second, we present an example scenario where the principal has to share its objective with the agent to induce the desired behavior. Third, we provide sufficient conditions on the structure of the MDP and the length of the agent's decision horizon under which there exists an optimal incentive design that allows the principal to hide its objective from the agent. 

\noindent \textbf{Related work.} The problem of obtaining desired agent behavior through a sequence of incentives has been extensively studied in the literature. In \cite{Zhang_value} and \cite{Zhang_general}, the authors present methods to design incentive sequences with \textit{limited} resources that maximizes the value of the principal's objective function. They employ techniques from inverse reinforcement learning literature and prove NP-hardness of the considered design problem \cite{Zhang_value}. The work \cite{Zhang} provides a polynomial-time algorithm to synthesize minimum incentives for inducing a \textit{specific} agent policy. Reference \cite{Chen_adaptive} considers a bandit model and presents methods to induce desired agent actions under different constraints on the incentives. Although it is quite different from the problem considered here, the design of feasible incentives that aligns the objectives of an agent and a principal is discussed in \cite{Control} from a control theoretic perspective. Unlike the references mentioned above, in this paper, we consider the problem of designing minimum incentives that maximizes the value of the principal's objective function expressed as a temporal logic formula. We also note that establishing the complexity of the design problem considered in this paper is mentioned as an open problem in \cite{Zhang_value}. 

\section{Preliminaries}

 For a set $S$, we denote its power set and cardinality by $2^S$ and $\lvert S \rvert$, respectively. Additionally, $\mathbb{N}$$=$$\{1,2,\ldots\}$, $\mathbb{N}_0$$=$$\{0,1,2,\ldots\}$ and $\mathbb{R}_{\geq 0}$$=$$[0,\infty)$.

\subsection{Markov Decision Processes}
{\setlength{\parindent}{0cm}
\begin{definition}
A \textit{Markov decision process} (MDP) is a tuple $\mathcal{M}$$=$$(S, s_0, \mathcal{A}, \mathcal{P},\mathcal{AP},\mathcal{L},\mathcal{R})$ where $S$ is a finite set of states, $s_0$$\in$$S$ is an initial state, $\mathcal{A}$ is a finite set of actions, $ \mathcal{P}$$:$$S$$\times$$ \mathcal{A}$$\times$$S$$\rightarrow$$[0,1]$ is a transition function such that $\sum_{s'\in S}\mathcal{P}(s,a,s')$$=$$1$ for all $s$$\in$$S$ and $a$$\in$$\mathcal{A}(s)$ where $\mathcal{A}(s)$ denote the available actions in $s$, $\mathcal{AP}$ is a set of atomic propositions, $\mathcal{L}$$:$$ S$$\rightarrow$$ 2^{\mathcal{AP}}$ is a function that labels each state with a subset of atomic propositions, and $\mathcal{R}$$:$$S$$\times$$\mathcal{A}$$\rightarrow$$\mathbb{R}$ is a reward function.
\end{definition}}
We denote the transition probability $ \mathcal{P}(s,a,s')$ by $ \mathcal{P}_{s,a,s'}$. 

{\setlength{\parindent}{0cm}
\noindent \begin{definition}
For an MDP $\mathcal{M}$, a \textit{decision rule} $d $$:$$ S$$\times$$\mathcal{A} $$\rightarrow$$[0,1]$ is a function such that $\sum_{a\in\mathcal{A}(s)}d(s,a)$$=$$1$ for all $s$$\in$$S$. A decision rule $d$ is said to be \textit{deterministic} if for all $s$$\in$$S$ there exists $a$$\in$$\mathcal{A}(s)$ such that $d(s,a)$$=$$1$, and \textit{randomized} otherwise. For an MDP $\mathcal{M}$, we denote the set of all (deterministic) decision rules by ($\mathcal{D}^D(\mathcal{M})$) $\mathcal{D}(\mathcal{M})$.
\end{definition}}

For an MDP $\mathcal{M}$, a decision-maker, i.e., an agent, chooses a decision rule $d$$\in$$\mathcal{D}(\mathcal{M})$ at each \textit{stage}. 

{\setlength{\parindent}{0cm}
\noindent \begin{definition}
An $N$-stage \textit{policy} for an MDP $\mathcal{M}$ is a sequence $\pi$$=$$(d_1, d_2, \ldots, d_N)$ where $N$$\leq$$\infty$ and $d_t$$\in$$\mathcal{D}(\mathcal{M})$ for all $t$$\leq$$N$.  A \textit{stationary} policy is a policy such that $d_t$$=$$d_1$ for all $t$$\leq$$N$. A policy is said to be \textit{deterministic} if $d_t$$\in$$\mathcal{D}^D(\mathcal{M})$ for all $t$, and \textit{randomized} otherwise.  For an MDP $\mathcal{M}$, we denote the set of all $N$-stage policies by $\Pi_N(\mathcal{M})$. For notational simplicity, we denote the set of $\infty$-stage policies by $\Pi(\mathcal{M})$.
\end{definition}}

For an MDP $\mathcal{M}$ and a policy $\pi$$\in$$\Pi(\mathcal{M})$, let $\mu^{\pi}_t(s,a)$ be the joint probability of being in state $s$$\in$$S$ and taking the action $a$$\in$$\mathcal{A}(s)$ at stage $t$, which is uniquely determined through the recursive formula
\begin{align}
    \mu_{t+1}^{\pi}(s',a')=\sum_{s\in S}\sum_{a\in \mathcal{A}(s)}\mathcal{P}_{s,a,s'} \mu_t^{\pi}(s,a)d_{t+1}(s',a')
\end{align}
where $\mu^{\pi}_1(s,a)$$=$$d_1(s,a)\mu_0(s)$ and $\mu_0$$:$$S$$\rightarrow$$\{0,1\}$ is a function such that $\mu_0(s_0)$$=$$1$ and $\mu_0(s)$$=$$0$ for all $s$$\in$$S\backslash \{s_0\}$.
{\setlength{\parindent}{0cm}
\noindent \begin{definition}
For an MDP $\mathcal{M}$ and a policy $\pi$$\in$$\Pi(\mathcal{M})$, the \textit{expected residence time} in a state-action pair $(s,a)$ is 
\begin{align}
\xi^{\pi}(s,a):=\sum_{t=1}^{\infty}\mu^{\pi}_t(s,a).
\end{align}
\end{definition}}

An infinite sequence $\varrho^{\pi}$$=$$s_0s_1s_2\ldots$ of states generated in $\mathcal{M}$ under a policy $\pi$$\in$$\Pi(\mathcal{M})$, which starts from the initial state $s_0$ and satisfies $\sum_{a_t\in \mathcal{A}(s_t)}d_k(s_t,a_t)\mathcal{P}_{s_t,a_t,s_{t+1}}$$>$$0$ for all $t$$\geq$$0$, is called a \textit{path}. Any finite prefix of $\varrho^{\pi}$ a finite path fragment. We define the set of all paths and finite path fragments in $\mathcal{M}$ under the policy $\pi$ by $Paths^{\pi}(\mathcal{M})$ and $Paths_{fin}^{\pi}(\mathcal{M})$, respectively. We use the standard probability measure over the outcome set $Paths^{\pi}(\mathcal{M})$ \cite{Model_checking}. 

{\setlength{\parindent}{0cm}
\noindent \begin{definition}
An \textit{incentive design} for an MDP $\mathcal{M}$
is a sequence $\Gamma$$=$$(\gamma_1,\gamma_2,\ldots)$ where $\gamma_t$$:$$S$$\times $$\mathcal{A}$$\rightarrow$$\mathbb{R}_{\geq 0}$. A \textit{stationary} incentive design is a design such that $\gamma_t$$=$$\gamma_1$ for all $t$$\in$$\mathbb{N}$. For an MDP $\mathcal{M}$, we denote the set of all incentive designs by $\Theta(\mathcal{M})$.
\end{definition}}
\subsection{Linear temporal logic}
We consider syntactically co-safe linear temporal logic (scLTL) formulae to specify tasks and refer the reader to \cite{Model_checking,Belta_book} for the syntax and semantics of scLTL. 

An scLTL formula is built up from a set $\mathcal{AP}$ of atomic propositions, logical connectives such as conjunction ($\land$) and negation ($\lnot$), and temporal modal operators such as until ($\mathcal{U}$) and eventually ($\lozenge$). An infinite sequence of subsets of $\mathcal{AP}$ defines an infinite \textit{word}, and an scLTL formula is interpreted over infinite words on $2^{\mathcal{AP}}$. We denote by $w$$\models$$\varphi$ that a word $w$$=$$w_0w_1w_2\ldots$ satisfies an scLTL formula $\varphi$.

For an MDP $\mathcal{M}$ under a policy $\pi$, a path $\varrho^{\pi}$$=$$s_0s_1\ldots$ generates a word $w$$=$$w_0w_1\ldots$ where $w_k$$=$$\mathcal{L}(s_k)$ for all $k$$\geq$$0$. With a slight abuse of notation, we use $\mathcal{L}(\varrho^{\pi})$ to denote the word generated by $\varrho^{\pi}$. For an scLTL formula $\varphi$, the set $\{\varrho^{\pi}$$\in$$ Paths^{\pi}(\mathcal{M})$$:$$ \mathcal{L}(\varrho^{\pi})$$\models$$\varphi\}$ is measurable \cite{Model_checking}. Hence, we define 
\begin{equation*}
\begin{aligned}
\text{Pr}_{\mathcal{M}}^{\pi}(s_0\models \varphi):=\text{Pr}_{\mathcal{M}}^{\pi}\{\varrho^{\pi}\in Paths^{\pi}(\mathcal{M}) : \mathcal{L}(\varrho^{\pi})\models \varphi\}
\end{aligned}
\end{equation*}
as the probability of satisfying the scLTL formula $\varphi$ for an MDP $\mathcal{M}$ under the policy $\pi$$\in$$\Pi(\mathcal{M})$.


\section{Problem Statement}\label{prob_state}
We consider an \textit{agent} whose sequential decision-making process is modeled as an MDP $\mathcal{M}$, and a \textit{principal} that provides the agent a sequence of incentives $\Gamma$$\in$$\Theta(\mathcal{M})$.

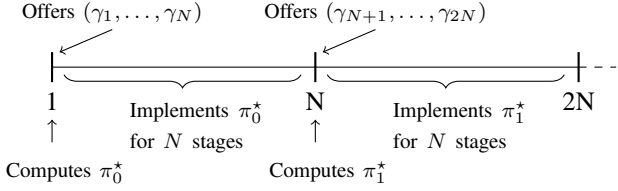
\begin{figure}[t]
\begin{tikzpicture}
\draw (0,0) -- (7,0);
\draw (7,0)[dashed] -- (7.5,0);
\draw[thick] (0,0.2) -- (0,-0.2) node[below] {1};
\draw[->] (0.8,0.5) to (0.1,0.2) node[xshift=20pt, yshift=15pt] {\footnotesize{Offers $(\gamma_1,\ldots,\gamma_N)$}};
\draw[->] (4.3,0.5) to (3.6,0.2) node[xshift=20pt, yshift=15pt] {\footnotesize{Offers $(\gamma_{N+1},\ldots,\gamma_{2N})$}};
\draw[->] (0,-1) to (0,-0.7) node[xshift=5pt, yshift=-20pt] {\footnotesize{Computes $\pi^{\star}_0$}};
\draw[thick] (3.5,0.2) -- (3.5,-0.2) node[below] {N};
\draw[->] (3.5,-1) to (3.5,-0.7) node[xshift=5pt, yshift=-20pt] {\footnotesize{Computes $\pi^{\star}_1$}};
\draw[thick] (7,0.2) -- (7,-0.2) node[below] {2N};
\draw [decorate,decoration={brace,mirror,amplitude=5pt},xshift=-4pt,yshift=0pt]
(0.3,-0.1) -- (3.45,-0.1) node [black,midway,xshift=0.3cm,yshift=-0.7cm, text width=2cm] 
{\footnotesize Implements $\pi^{\star}_0$ \\ for $N$ stages};
\draw [decorate,decoration={brace,mirror,amplitude=5pt},xshift=-4pt,yshift=0pt]
(3.8,-0.1) -- (6.95,-0.1) node [black,midway,xshift=0.3cm,yshift=-0.7cm, text width=2cm] 
{\footnotesize Implements $\pi^{\star}_1$ \\ for $N$ stages};
\end{tikzpicture}
\caption{An illustration of the incentive implementation and the agent's decision-making process. The principal offers incentives for the next $N$ stages. After receiving the incentive offers, the agent computes and implements its optimal decisions for the next $N$ stages.}
\label{graphical_rep}\vspace{-0.3cm}
\end{figure}

The agent's objective is to maximize its expected total reward after $N$ stages. However, since the incentive sequence offered by the principal might be non-stationary, the agent computes an $N$-stage policy every $N$ stages. A graphical illustration of the agent's planning method is shown in Fig. \ref{graphical_rep}. Formally, let $N$$\in$$\mathbb{N}$ be a constant, and $\mathcal{R}(S_t,A_t)$ and $\gamma_{t}(S_t,A_t)$ be the random reward and incentive received in stage $t$$\leq$$N$. Additionally, let $J$$:=$$(J_0,J_1,\ldots)$ be a sequence of objective functions where $J_k$$:$$\Pi_N(\mathcal{M})$$\times$$\Theta(\mathcal{M})$$\rightarrow$$\mathbb{R}^{\lvert S\rvert}$ is such that 
\begin{align*}
J_k(\pi,\Gamma)(s):=\mathbb{E}_{s}^{\pi}\Big[\sum_{t=1}^{N} (\mathcal{R}(S_t,A_t)+\gamma_{kN+t}(S_t,A_t))\Big]
\end{align*}
for all $s$$\in$$S$ where the expectation is taken over the finite path fragments that are generated by the policy $\pi$$\in$$\Pi_N(\mathcal{M})$ and start from the state $s$. 
Then, for a given incentive design $\Gamma$$\in$$\Theta(\mathcal{M})$, the agent's optimal $\infty$-stage policy is given by $\pi^{\star}$$:=$$(\pi_0^{\star}, \pi_1^{\star}, \ldots)$ where $\pi^{\star}_k$ is such that
\begin{align}
\pi_k^{\star}\in\arg\max_{\pi\in \Pi_N(\mathcal{M})} J_k(\pi,\Gamma)(s)
\end{align}
for all $s$$\in$$S$ and $k$$\in$$\mathbb{N}_0$. Note that the agent's policy $\pi^{\star}_k$ maximizes the total reward starting from any $s$$\in$$S$.

The principal's objective is to design an incentive sequence such that the agent's optimal policy under the provided incentives satisfies an scLTL formula $\varphi$ with maximum probability. 

The problem that we consider is the synthesis of an incentive design that minimizes the cost to the principal while realizing its objective. We make the following assumptions:

\begin{enumerate}[(i)]
    \item Agent's reward function $\mathcal{R}$ is known by the principal.
    \item Agent's decision horizon $N$ is known by the principal.
    \item The principal pays the offered incentives if and only if the agent takes the incentivized action.
\end{enumerate}

Then, the optimization problem that we are interested in to solve is the following:
\begin{subequations}\label{problem}
\begin{align}\label{main_prob_1}
    \min_{\Gamma\in \Theta(\mathcal{M})}&\ \  \mathbb{E}_{s_0}^{\pi^{\star}}\Big[\sum_{t=1}^{\infty}\gamma_t(s,a)\Big]\\ \label{sup}
     \text{subject to:}& \ \ \pi^{\star}=(\pi^{\star}_0,\pi^{\star}_1,\ldots)\\
    & \pi_k^{\star}\in\arg\max_{\pi\in \Pi_N(\mathcal{M})}J_k(\pi,\Gamma)(s) \ \forall s\in S ,\forall k\in \mathbb{N}_0\\
    &\ \ \text{Pr}_{\mathcal{M}}^{\pi^{\star}}(s_0\models\label{maxreach} \varphi)=\max_{\pi\in\Pi(\mathcal{M})}\text{Pr}_{\mathcal{M}}^{\pi}(s_0\models \varphi)
\end{align}
\end{subequations}
where $\Gamma$$=$$(\gamma_1,\gamma_2,\ldots)$.


\section{The Design of Incentive Sequences}\label{thecostofcontrol}
In this section, we provide a method to synthesize an incentive design that solves the problem \eqref{main_prob_1}-\eqref{maxreach}. For simplicity, we restrict our attention to reachability specifications, i.e., $\varphi$$=$$\lozenge p$ where $p$$\in$$\mathcal{AP}$. The incentive design for general scLTL specifications is discussed in Section \ref{Generalization}.

We first partition the states into three disjoint sets as follows. Let $B$$\subseteq$$S$ be the set of all states such that $\{p\}$$\subseteq$$\mathcal{L}(s)$, i.e., the set of states that the principal wants the agent to reach, and $S_0$$\subseteq$$S$ be the set of states that have zero probability of reaching the states in $B$ under any policy. More precisely, $s$$\in$$S_0$ if $\text{Pr}_{\mathcal{M}}^{\pi}(s\models \lozenge p)$$=$$0$ for all $\pi$$\in$$\Pi(\mathcal{M})$. Finally, we let $S_r$$=$$S\backslash B\cup S_0$ be the set of all states that are not in $B$ and have nonzero probability of reaching a state in $B$ under some policy. These sets can be found in time polynomial in the size of the MDP using graph search algorithms \cite{Model_checking}.

The agent's initial state $s_0$$\in$$S$ can belong to either $B$, $S_0$ or $S_r$. However, we only consider the case $s_0$$\in$$S_r$ since otherwise the optimal incentive design is trivially $\gamma_t(s,a)$$=$$0$ for all $t$$\in$$\mathbb{N}$.


\subsection{The cost of control}\label{caseN2}
Recall that the agent's first objective function $J_0$$:$$\Pi_N(\mathcal{M})$$\times$$\Theta(\mathcal{M})$$\rightarrow$$\mathbb{R}^{\lvert S \rvert}$ is 
\begin{align*}
J_0(\pi,\Gamma)(s)=\mathbb{E}_{s}^{\pi}\Big[\sum_{t=1}^{N} (\mathcal{R}(S_t,\mathcal{A}_t)+\gamma_{t}(S_t,\mathcal{A}_t))\Big]
\end{align*}
for all $s$$\in$$S$. Let $V_n$$:$$S$$\rightarrow$$\mathbb{R}$ be the agent's \textit{value function} at stage $n$ such that
\begin{align*}
    V_n(s):=\max_{\pi\in\Pi_N(\mathcal{M})}\mathbb{E}_{s}^{\pi}\Big[\sum_{t=n}^{N} (\mathcal{R}(S_t,\mathcal{A}_t)+\gamma_{t}(S_t,\mathcal{A}_t))\Big]
\end{align*}
for all $s$$\in$$S$, where the expectation is taken over the paths that occupy $s$ at stage $n$.  Then, we have the recursive formula
\begin{align*}
    V_n(s)=\max_{a\in \mathcal{A}(s)}\mathcal{R}(s,a)+\gamma_{n}(s,a)+\sum_{s'\in S}\mathcal{P}_{s,a,s'}V_{n+1}(s')
\end{align*}
for all $1$$\leq$$n$$\leq$$N$, where  $V_{N+1}(s)$$=$$0$ for all $s$$\in$$S$. Let $Q_n$$:$$S$$\times$$\mathcal{A}$$\rightarrow$$\mathbb{R}$ be the agent's \textit{$Q$-function} at stage $n$ such that 
\begin{align*}
    Q_n(s,a):=\mathcal{R}(s,a)+\gamma_{n}(s,a)+\sum_{s'\in S}\mathcal{P}_{s,a,s'}V_{n+1}(s').
\end{align*}
By the principle of optimality \cite{Puterman},\cite{Bertsekas}, the agent's optimal policy $\pi^{\star}_0$$=$$(d^{\star}_1,d^{\star}_2,\ldots,d^{\star}_N)$ is such that, for all $1$$\leq$$n$$\leq$$N$, $d_n(s, a')$$>$$0$ only if 
\begin{align*}
    a'\in\arg\max_{a\in \mathcal{A}(s)}Q_n(s,a).
\end{align*}
We recursively define
\begin{align}\label{q1}
\overline{Q}_n(s,a)&:=\mathcal{R}(s,a)+\sum_{s'\in S}\mathcal{P}_{s,a,s'}\overline{V}_{n+1}(s'),\\ \label{v1}
        \overline{V}_n(s)&:=\max_{a\in \mathcal{A}(s)}\overline{Q}_n(s,a),
\end{align}
for all $s$$\in$$S$ and $a$$\in$$\mathcal{A}(s)$.
For a given $\epsilon$$\geq$$0$, we finally define a real-valued function $\phi^{\epsilon}_t$$:$$S\times \mathcal{A}$$\rightarrow$$\mathbb{R}_{\geq0}$ such that
\begin{align*}
   \phi^{\epsilon}_n(s,a):=\begin{cases} \overline{V}_n(s)-\overline{Q}_n(s,a)+\epsilon & \text{if}\ s\in S_r, \ a\in\mathcal{A}(s)\\
   0 & \text{otherwise}.
    \end{cases}
\end{align*}

For an arbitrarily small $\epsilon$$>$$0$, the value of $\phi^{\epsilon}_n(s,a)$, referred as \textit{the cost of control} for the state-action pair $(s,a)$, is the minimum incentive that should be offered to the agent in order to make the action $a$$\in$$\mathcal{A}(s)$ uniquely optimal at stage $t$. It is worth noting that although the cost of control $\phi^{\epsilon}_n(s,a)$ depends on the stage number $n$, it is independent of the objective number, i.e., it is the same for all $J_k$. This is because the agent's reward function $\mathcal{R}$ is stationary, and therefore, $\overline{V}_n(s)$ and $\overline{Q}_n(s,a)$ do not change with the objective number $k$ as can be seen from \eqref{q1}-\eqref{v1}.

\subsection{An $\overline{\epsilon}$-optimal incentive design}\label{expandedMDPsection}

To synthesize the minimum incentive sequence, we should specify the actions to be incentivized by the principal at each state for each stage. To this aim, we modify the MDP $\mathcal{M}$ by considering the agent's decision horizon $N$ as another dimension in the state-space.

{\setlength{\parindent}{0cm}
\noindent \begin{definition}
For an MDP $\mathcal{M}$ and $T$$=$$\{1,2,\ldots,N\}$, the \textit{expanded MDP} is a tuple $\overline{\mathcal{M}}$$=$$(\overline{S},\overline{s}_0 , \mathcal{A}, \overline{\mathcal{P}},\mathcal{AP},\overline{\mathcal{L}},\mathcal{R})$ where
\begin{itemize}
    \item $\overline{S}$$=$$S\times T$,
    \item $\overline{s}_0$$=$$(s_0,1)$ is the initial state,
    \item $\overline{\mathcal{P}}$$:$$\overline{S} \times \mathcal{A} \times \overline{S}$$\rightarrow$$[0,1]$ is such that
    \begin{align*}
    &\overline{\mathcal{P}}_{(s,n),a,(s',n')}=\\
    &\quad\qquad  \begin{cases} \mathcal{P}_{s,a,s'} & \text{if} \ 1\leq n\leq N-1\  \text{and} \  n'=n+1\\ \mathcal{P}_{s,a,s'} & \text{if} \ n= N\  \text{and} \  n'=1\\
    0 & \text{otherwise},\end{cases}\end{align*}
    \item $\overline{\mathcal{L}}$$:$$\overline{S}$$\rightarrow$$2^{\mathcal{AP}}$ is such that $\overline{\mathcal{L}}((s,t))$$=$$\mathcal{L}(s)$ for all $s$$\in$$S$ and for all $t$$\in$$T$,
\end{itemize}
and $\mathcal{A}$, $\mathcal{AP}$ and $\mathcal{R}$ are as defined for $\mathcal{M}$.
\end{definition}}

We note that the transition function $\overline{\mathcal{P}}$ is defined such that the agent's initial state while computing the $k$-th $N$ stage policy is the state occupied by the agent at $kN$$+$$1$-st stage on the expanded MDP. 

Let $\overline{B}$$\cup$$\overline{S}_0$$\cup$$\overline{S}_r$ be the partition of the states of  $\overline{\mathcal{M}}$ such that if $s$$\in$$B$, then $(s,n)$$\in$$\overline{B}$ for all $n$$\in$$T$, and the sets $\overline{S}_0$ and $\overline{S}_r$ are defined similarly. Then, the principal's objective on $\overline{\mathcal{M}}$ is to induce an agent policy that reaches the set $\overline{B}$ with maximum probability. To synthesize an incentive design under which the optimal agent policy satisfies the desired property, we modify the expanded MDP $\overline{\mathcal{M}}$ by making its states $s$$\in$$\overline{B}$$\cup$$\overline{S}_0$ absorbing, and denote the resulting MDP by $\overline{\mathcal{M}}'$. Then, for a given $\epsilon$$\geq$$0$, we define the cost of control for a state-action pair on $\overline{\mathcal{M}}$ through the function
$\underline{\phi}^{\epsilon}$$:$$\overline{S}$$\times$$ \mathcal{A}$$\rightarrow$$\mathbb{R}_{\geq0}$ such that
\begin{align*}
    \underline{\phi}^{\epsilon}((s,n),a):=\begin{cases} \overline{V}_n(s)-\overline{Q}_n(s,a)+\epsilon & \text{if}\ s\in \overline{S}_r, a\in\mathcal{A}(s)\\
   0 & \text{otherwise}.
    \end{cases}
\end{align*}
Let $\Xi(\overline{\mathcal{M}}')$$\subseteq$$\Pi(\overline{\mathcal{M}}')$
be a subset of the set of $\infty$-stage policies such that $\pi'$$\in$$\Xi(\overline{\mathcal{M}}')$ if and only if 
\begin{align}
    \pi'\in\arg\max_{\pi\in\Pi(\overline{\mathcal{M}}')}\text{Pr}^{\pi}(\overline{s}_0\models \varphi),
\end{align}
and for $\epsilon$$\geq$$0$, $f_{\epsilon}$ $:$ $\Xi(\overline{\mathcal{M}}')$$\rightarrow$$\mathbb{R}$ be a function such that 
 \begin{align}\label{objfunction}
    f_{\epsilon}(\pi):=\mathbb{E}_{\overline{s}_0}^{\pi}\Big[\sum_{t=1}^{\infty}\underline{\phi}^{\epsilon}(S_t,\mathcal{A}_t)\Big].
\end{align}
Then, for an arbitrarily small $\overline{\epsilon}$$>$$0$, an $\overline{\epsilon}$-optimal incentive sequence can be designed in two steps as follows. 

\textit{Step 1:} Compute $\overline{V}_n(s)$ and $\overline{Q}_n(s,a)$ given in \eqref{q1}-\eqref{v1}, and construct the cost of control function $\underline{\phi}^{\epsilon}$. Then for the modified expanded MDP $\overline{\mathcal{M}}'$, compute a \textit{stationary deterministic} policy $\widetilde{\pi}$$=$$(\widetilde{d},\widetilde{d},\ldots)$ such that
\begin{align}
    \widetilde{\pi}\in\arg\min_{\pi\in\Xi(\overline{\mathcal{M}}')}f_{\epsilon}(\pi).
\end{align}

\textit{Step 2:} Let $\varrho^{\pi}$$\in$$Paths^{\pi}(\mathcal{M})$ be the path followed by the agent. At stage $kN$ where $N$ is the agent's decision horizon and $k$$\in$$\mathbb{N}_0$, provide the agent with the incentive sequence $\{\widetilde{\gamma}_1,\widetilde{\gamma}_2, \ldots, \widetilde{\gamma}_N\}$ such that 
\begin{itemize}
    \item if $\varrho^{\pi}[n]$$\not\in$$ B\cup S_0$ for all $n$$\leq$$kN$
\end{itemize}
\begin{align}\label{optincentives}
    \widetilde{\gamma}_n(s,a):=\begin{cases}
    \underline{\phi}^{\epsilon}((s,n),a) & \text{if} \  s\in S_r \ \text{and}  \ \widetilde{d}((s,n))(a)>0,  \\
    \epsilon & \text{if} \  s\not\in S_r \ \text{and} \ \widetilde{d}((s,n))(a)>0, \\
    0 & \text{otherwise},
    \end{cases}
\end{align}
\begin{itemize}
    \item $\widetilde{\gamma}_n(s,a)$$:=$$0$ otherwise.
\end{itemize}

Under the proposed incentive design \eqref{optincentives}, the agent's value function $V_n$ satisfies $V_n(s)$$=$$\overline{V}_n(s)$$+$$(N+1-n)\epsilon$ for all $s$$\in$$S$, $n$$\leq$$N$. Additionally, if $\varrho^{\pi}[n]$$\not\in$$ B\cup S_0$ for all $n$$\leq$$kN$, then for all $s$$\in$$S$, $\widetilde{d}((s,n))(a)$$>$$0$ implies that the agent's $Q$-function satisfies
\begin{align*}
        Q_n(s,a)&=\mathcal{R}(s,a)+\widetilde{\gamma}_{n}(s,a)+\sum_{s'\in S}\mathcal{P}_{s,a,s'}V_{n+1}(s')\\
        &=\widetilde{\gamma}_{n}(s,a)+\overline{Q}_n(s,a)+(N-n)\epsilon\\
        &=(N+1-n)\epsilon+\overline{V}_n(s)\\
        &>(N-n)\epsilon+\overline{V}_n(s)=\max_{a'\in\mathcal{A}(s)\backslash\{a\}}Q_n(s,a').
\end{align*}
Consequently, the agent is guaranteed to take the incentivitized actions at each stage until reaching the set $B$$\cup$$S_0$.

We now show $\overline{\epsilon}$-optimality of the proposed incentive design. Note that an optimal incentive design, i.e., $\overline{\epsilon}$$=$$0$, does not exist since choosing $\epsilon$$=$$0$ in the cost of control function $\phi^{\epsilon}_n$ may not make the incentivized action uniquely optimal for the agent. As a result, the principal may not be able to control the agent's actions by offering such incentives.

We need the following technical lemma to state the main result. 
{\setlength{\parindent}{0cm}
\begin{lemma}\label{finiteness_lemma} There exists a policy $\widetilde{\pi}$$\in$$\arg\min_{\pi\in\Xi(\overline{\mathcal{M}}')}f_{0}(\pi)$ such that $\xi^{\widetilde{\pi}}(s,a)$$<$$\infty$ for all $s$$\in$$\overline{S}_r$ and $a$$\in$$\mathcal{A}(s)$.
\end{lemma}}

\noindent\textbf{Proof (Sketch):} The problem of synthesizing a policy $\widetilde{\pi}$ such that $\widetilde{\pi}$$\in$$\arg\min_{\pi\in\Xi(\overline{\mathcal{M}}')}f_{0}(\pi)$ can be recast as a stochastic shortest path (SSP) problem with dead ends and zero-cost loops. Specifically, the dead ends are the states $\overline{S}_0$ and zero-cost loops are formed by states $\overline{S}_r$. The existence of stationary policies for such SSP problems can be established by slightly modifying the statement of Theorem 1 in \cite{SSP1}. Since any stationary policy $\pi$$\in$$\Xi(\overline{\mathcal{M}}')$ is guaranteed to reach the set $\overline{B}\cup \overline{S}_0$ with probability 1 within finite number of stages, the result follows. $\Box$

{\setlength{\parindent}{0cm}
\begin{thm} \label{epsopttheorem} For any given $\overline{\epsilon}$$>$$0$, there exists $\epsilon$$>$$0$ such that
\begin{align*}
    \min_{\pi\in\Xi(\overline{\mathcal{M}}')} f_{\epsilon}(\pi)\leq \min_{\pi\in\Xi(\overline{\mathcal{M}}')} f_0(\pi)+\overline{\epsilon}.
\end{align*}
\end{thm}}

\noindent\textbf{Proof:} For any policy $\pi$$\in$$\Xi(\overline{\mathcal{M}}')$ such that $\xi^{\pi}(s,a)$$<$$\infty$ for all $s$$\in$$\overline{S}_r$ and $a$$\in$$\mathcal{A}(s)$, we have
\begin{align}\label{eqeq}
    & f_{\epsilon}(\pi)=  f_0(\pi)+\sum_{s\in \overline{S}_r}\sum_{a\in\mathcal{A}(s)}\xi^{\pi}(s,a)\epsilon.
\end{align}
Now, for a given $\overline{\epsilon}$$>$$0$, we evaluate both sides of the above equation at $\overline{\pi}$$\in$$\arg\min_{\pi\in\Xi(\overline{\mathcal{M}}')} f_0(\pi)$, which satisfies the condition $\xi^{\overline{\pi}}(s,a)$$<$$\infty$ due to Lemma \ref{finiteness_lemma}. Choosing
\begin{align*}
\epsilon=\frac{\overline{\epsilon}}{\sum_{s\in S_r}\sum_{a\in\mathcal{A}(s)}\xi^{\overline{\pi}}(s,a)}>0
\end{align*}
and taking the minimum of the left hand side of \eqref{eqeq} over the set $\Xi(\mathcal{M})$, we conclude the result. $\Box$

We conclude this section by noticing a remarkable property of the proposed incentive design. Specifically, to implement the proposed design \eqref{optincentives}, the principal should use only a simple switch mode which offers the same incentives until the agent reaches the set $B$$\cup$$S_0$ and shifts all incentives to zero after the agent either satisfies the principal's objective or fails to satisfy it. 


\section{Computation of an Optimal Incentive Design}\label{computational_section}
In the previous section, we developed a method to synthesize an $\overline{\epsilon}$-optimal incentive design which require us to solve a constrained cost minimization problem given in \eqref{objfunction}. Specifically, to solve the incentive design problem \eqref{main_prob_1}-\eqref{maxreach}, one should synthesize a \textit{stationary deterministic} policy $\widetilde{\pi}$ such that
\begin{align}\label{main_problem}
    \widetilde{\pi}\in\arg\min_{\pi\in\Xi(\overline{\mathcal{M}}')}\mathbb{E}_{s_0}^{\pi}\Big[\sum_{t=1}^{\infty}\underline{\phi}^{\epsilon}(S_t,\mathcal{A}_t)\Big]
\end{align}

In this section, we develop a method to solve the above optimization problem. For the ease of notation, we consider an scLTL formula of the form $\varphi$$=$$\lozenge p$. The incentive design for general scLTL formulae is discussed in Section \ref{Generalization}.

\subsection{Construction of the feasible policy space}
To solve the problem \eqref{main_problem}, we first represent the set $\Xi(\overline{\mathcal{M}}')$ of feasible policies as a set of policies that maximizes the expected total reward with respect to a specific reward function. 

For a given MDP $\overline{\mathcal{M}}$, we partition the set of states into three disjoint sets $\overline{B}$, $\overline{S}_0$, and $\overline{S}_r$ as explained in Section \ref{thecostofcontrol}, and make the states $s$$\in$$\overline{B}$$\cup$$\overline{S}_0$ absorbing to form the modified MDP $\overline{\mathcal{M}}'$. For the modified MDP, we define a reward function $r$$:$$\overline{S}$$\times$$\mathcal{A}$$\rightarrow$$\mathbb{R}_{\geq 0}$ such that
\begin{align*}
    r(s,a)=\begin{cases}
    \sum_{s'\in \overline{B}}\overline{\mathcal{P}}_{s,a,s'} & \text{if} \ \ s\in \overline{S}_r\\
    0 & \text{otherwise.}
    \end{cases}
\end{align*}
By making use of the known results, e.g., Theorem 10.100 in \cite{Model_checking}, it can be easily shown that for any $s$$\in$$\overline{S}$ and $\pi$$\in$$\Pi(\mathcal{M}')$, 
\begin{align*}
    \mathbb{E}_{s}^{\pi}\Big[\sum_{t=1}^{\infty}r(S_t,\mathcal{A}_t)\Big]=\text{Pr}^{\pi}(s\models \varphi)
\end{align*}
where $\varphi$$=$$\lozenge p$, $p$$\in$$\mathcal{AP}$, and $\{p\}$$\subseteq$$\mathcal{L}(s')$ if and only if $s'$$\in$$\overline{B}$. Let $x_s^{\star}$$:=$$\max_{\pi\in\Pi(\overline{\mathcal{M}}')}\text{Pr}^{\pi}(s\models \varphi)$. Then, the problem \eqref{main_problem} can be rewritten as
\begin{subequations}\label{main_problem2}
\begin{align}\label{probprob}
     \min_{\pi\in\Pi(\overline{\mathcal{M}}')}&\ \ \mathbb{E}_{s_0}^{\pi}\Big[\sum_{t=1}^{\infty}\underline{\phi}^{\epsilon}(S_t,\mathcal{A}_t)\Big]\\ \label{constraint}
     \text{subject to:}&\ \  \mathbb{E}_{s_0}^{\pi}\Big[\sum_{t=1}^{\infty}r(S_t,\mathcal{A}_t)\Big]=x_{s_0}^{\star}.
\end{align}
\end{subequations}

\subsection{Synthesis of an optimal stationary deterministic policy}
Using Lemma \ref{finiteness_lemma}, one can formulate the problem \eqref{probprob}-\eqref{constraint} as a linear optimization problem and synthesize an optimal stationary policy. First, we compute the maximum probability of satisfying the specification $\varphi$, i.e., $x^{\star}_{s_0}$$=$$\max_{\pi\in\Pi(\overline{\mathcal{M}}')}\text{Pr}^{\pi}(s_0\models \varphi)$, by solving a linear program (LP) \cite{Model_checking} (see Chapter 10). Then we solve the following LP
\begin{subequations}\vspace{-0.2cm}\label{opt_1}
\begin{flalign}\label{probobj}
&\underset{\lambda(s,a)}{\text{minimize}}\qquad\sum_{s\in \overline{S}_r}\sum_{a\in \mathcal{A}} \lambda(s,a)\underline{\phi}^{\epsilon}(s,a)\\ \label{cons_1}
&\text{subject to:} \qquad\sum_{s\in \overline{S}_r}\sum_{a\in \mathcal{A}} \lambda(s,a)r(s,a)=x^{\star}_{s_0}\\ \label{cons_2}
&\forall s\in \overline{S}_r, \ \sum_{a\in \mathcal{A}(s)} \lambda(s,a)-\sum_{s'\in \overline{S}_r}\sum_{a\in \mathcal{A}(s)} \overline{\mathcal{P}}_{s',a,s}\lambda(s',a)=\alpha(s)  \\ \label{cons_last}
&\forall s\in \overline{S}_r, \ a\in\mathcal{A}(s),\  \lambda(s,a)\geq 0  \ 
\end{flalign}
\end{subequations}
where $\alpha$$:$$\overline{S}$$\rightarrow$$\{0,1\}$ is a function such that $\alpha(s_0)$$=$$1$ and $\alpha(s)$$=$$0$ for all $s$$\in$$\overline{S}\backslash\{\overline{s}_0\}$. The variable $\lambda(s,a)$ denotes the expected residence time in the state-action pair $(s,a)$ \cite{Marta,Puterman}. The constraint \eqref{cons_1} ensures that the probability of satisfying the specification $\varphi$ is maximized, and the constraints \eqref{cons_2} represent the balance between the \say{inflow} to and \say{outflow} from states.

For each $s$$\in$$\overline{S}_r$ and $a$$\in$$\mathcal{A}(s)$, let $\lambda^{\star}(s,a)$ be optimal decision variables in \eqref{probobj}-\eqref{cons_last}. An optimal stationary policy $\pi^{\star}$$=$$\{d^{\star},d^{\star},\ldots\}$ that solves the problem \eqref{probprob}-\eqref{constraint} is then given by 
\begin{align}\vspace{-0.2cm}\label{classic_randomized}
    d^{\star}(s,a):=\begin{cases} 
    \frac{\lambda^{\star}(s,a)}{\sum_{a\in\mathcal{A}(s)}\lambda^{\star}(s,a)} & \text{if} \  \sum_{a\in\mathcal{A}(s)}\lambda^{\star}(s,a)>0\\
    \text{arbitrary} & \text{otherwise}\end{cases}
    \raisetag{13pt}
\end{align}
for $s$$\in$$\overline{S}_r$, and $d^{\star}(s,a)$$=$$1$ for an arbitrary $a$$\in$$\mathcal{A}(s)$ for $s$$\not\in$$\overline{S}_r$.

We note that a policy constructed through \eqref{classic_randomized} is randomized in general. One can argue that choosing one of the actions $a$$\in$$\mathcal{A}(s)$ such that $ d^{\star}(s,a)$$>$$0$ deterministically yields an optimal stationary deterministic policy. However, the following example illustrates that such an approach may result in an infeasible policy for the problem \eqref{probobj}-\eqref{cons_last}.

\begin{figure}[b!]\centering\vspace{-0.5cm}
\scalebox{0.9}{\begin{tikzpicture}[->, >=stealth', auto, semithick, node distance=2cm]

    \tikzstyle{every state}=[fill=white,draw=black,thick,text=black,scale=0.7]

    \node[state,initial,initial text=] (s_0) {$s_0$};
    \node[state] (s_1) [right=10mm of s_0]  {$s_1$};
     \node[state] (s_2) [right=10mm of s_1]  {$s_2$};

\path
(s_0)  edge[out=290, in=250]  node[below]{$a_1, 0$}     (s_1)
(s_1)  edge[out=110, in=70]      node[above]{$a_1, 0$}     (s_0)
(s_1)  edge   node{$a_2, 1$}     (s_2)
(s_2)  edge  [loop right=10]    node{$a_1,0$}     (s_2);
\end{tikzpicture}}
\caption{ An MDP example for which arbitrarily choosing one of the optimal actions and taking it deterministically yields an infeasible policy. }\label{randomized_example}
\end{figure}
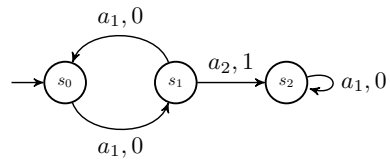

\noindent \textbf{Example 1:} Consider the MDP given in Fig. \ref{randomized_example}, where the cost of control $\underline{\phi}^{\epsilon}$ is such that $\underline{\phi}^{\epsilon}(s_1,a_2)$$=$$1$ and $\underline{\phi}^{\epsilon}(s,a)$$=$$0$ otherwise. Suppose that the specification is $\varphi$$=$$\lozenge s_2$, i.e., $r(s_1,a_2)$$=$$1$ and $r(s,a)$$=$$0$ otherwise. For the LP \eqref{probobj}-\eqref{cons_last}, a set of optimal decision variables is given by $\lambda^{\star}(s_0,a_1)$$=$$2$, $\lambda^{\star}(s_1,a_1)$$=$$1$, and $\lambda^{\star}(s_1,a_2)$$=$$1$. Therefore, an optimal policy synthesized through \eqref{classic_randomized} is $d^{\star}(s_0,a_1)$$=$$1$, $d^{\star}(s_1,a_1)$$=$$1/2$, and $d^{\star}(s_1,a_2)$$=$$1/2$. Clearly, if we consider a deterministic policy such that $d(s_1,a_1)$$=$$1$, the probability of satisfying the specification $\varphi$ under this policy is zero. Hence, choosing an arbitrary action $a$$\in$$\mathcal{A}(s)$ such that $ d^{\star}(s,a)$$>$$0$ deterministically violates the constraint and yields an infeasible policy.$\triangleleft$

As Example 1 illustrates, a structured approach is required to synthesize an optimal deterministic policy from the solution of the LP \eqref{probobj}-\eqref{cons_last}. 
 Let $\upsilon^{\star}$ be the optimal value of the LP in \eqref{probobj}-\eqref{cons_last}. To synthesize an optimal deterministic policy, we first solve the following LP, 
\begin{subequations}\label{opt_2}
\begin{align}\label{obj22}
&\underset{\lambda(s,a)}{\text{minimize}}\qquad\sum_{s\in \overline{S}_r}\sum_{a\in \mathcal{A}} \lambda(s,a)\\ \label{cons_2_1}
&\text{subject to:} \quad\  \sum_{s\in \overline{S}_r}\sum_{a\in \mathcal{A}} \lambda(s,a)r(s,a)=x^{\star}_{s_0}\\ \label{cons_2_2}
&\qquad \qquad \ \ \ \ \sum_{s\in \overline{S}_r}\sum_{a\in \mathcal{A}} \lambda(s,a)\underline{\phi}^{\epsilon}(s,a)=\upsilon^{\star}\\
&\forall s\in \overline{S}_r, \sum_{a\in \mathcal{A}(s)} \lambda(s,a)-\sum_{s'\in \overline{S}_r}\sum_{a\in \mathcal{A}(s)} \overline{\mathcal{P}}_{s',a,s}\lambda(s',a)=\alpha(s)\\  \label{const_last2}
&\forall s\in \overline{S}_r, \ \ a\in\mathcal{A}(s), \ \lambda(s,a)\geq 0 .
\end{align}
\end{subequations}
From the optimal decision variables $\lambda^{\star}(s,a)$ of  \eqref{obj22}-\eqref{const_last2}, an optimal policy $\pi^{\star}$$=$$\{d^{\star},d^{\star},\ldots\}$ can be generated as follows. Let $\mathcal{A}^{\star}(s)$$:=$$\{a$$\in$$\mathcal{A}(s) : \lambda^{\star}(s,a)$$>$$0 \}$. If $\mathcal{A}^{\star}(s)$$\neq$$\emptyset$, we choose $d^{\star}(s,a)$$=$$1$ for an arbitrary $a$$\in$$\mathcal{A}^{\star}(s)$, and if $\mathcal{A}^{\star}(s)$$=$$\emptyset$, we choose $d^{\star}(s,a)$$=$$1$ for an arbitrary $a$$\in$$\mathcal{A}(s)$. 

{\setlength{\parindent}{0cm}
\begin{prop}\label{comp_prop}
A stationary deterministic policy generated from the optimal decision variables $\lambda^{\star}(s,a)$ of  \eqref{obj22}-\eqref{const_last2} is a solution to the problem \eqref{probprob}-\eqref{constraint}.
\end{prop}}

A proof of Proposition \ref{comp_prop} can be found in Appendix \ref{proof_appendix}. Intuitively, the LP in \eqref{obj22}-\eqref{const_last2} computes the minimum expected time to reach the set $\overline{B}$ with probability $x^{\star}_{s_0}$ with the cost of $\upsilon^{\star}$. Therefore, if $\lambda^{\star}(s,a)$$>$$0$, by taking the action $a$$\in$$\mathcal{A}(s)$, the agent has to \say{get closer} to the set $\overline{B}$ with nonzero probability. Otherwise, the minimum expected time to reach the set $\overline{B}$ would be strictly decreased. Consequently, by choosing an arbitrary action $a$$\in$$\mathcal{A}^{\star}(s)$, the agent is guaranteed to reach the set $\overline{B}$ with the desired probability.

\section{Incentive Design for General scLTL specifications}\label{Generalization}

In previous sections, we have developed methods to synthesize $\overline{\epsilon}$-optimal incentive designs for reachability specifications $\varphi$$=$$\lozenge p$. For such specifications, the principal induces the desired agent behavior by sharing only the incentive sequences with the agent. In other words, the principal does not have to inform the agent explicitly about the specification. In this section, we show that for general scLTL formulae, the problem \eqref{main_prob_1}-\eqref{maxreach} may not have a feasible solution, in which case the principal must share its objective with the agent to induce the desired behavior. 

To solve the problem \eqref{main_prob_1}-\eqref{maxreach} for general scLTL formulae, one needs to utilize the techniques from automata theory \cite{Model_checking}. In particular, we use the fact that for any scLTL formula $\varphi$ built up from $\mathcal{AP}$, we can construct a \textit{deterministic finite automata} (DFA) $A_{\varphi}$$=$$(\mathcal{Q}, q_0, 2^{\mathcal{AP}},\delta_{\varphi},\mathcal{F})$ where $\mathcal{Q}$ is a finite set of memory states, $2^{\mathcal{AP}}$ is the alphabet, $\delta_{\varphi}$$:$$\mathcal{Q}\times 2^{\mathcal{AP}}$$\rightarrow$$\mathcal{Q}$ is a transition function and $\mathcal{F}$$\subseteq$$\mathcal{Q}$ is the set of accepting states \cite{Belta_book}.
Then, after forming the expanded MDP $\overline{\mathcal{M}}$ for a given MDP $\mathcal{M}$ and a decision horizon $N$ as explained in Section \ref{expandedMDPsection}, one can construct the product MDP which is defined as follows.
{\setlength{\parindent}{0cm}
 \noindent \begin{definition}
Let $\overline{\mathcal{M}}$$=$$(\overline{S}, \overline{s}_0, \mathcal{A}, \overline{\mathcal{P}}, \mathcal{AP}, \overline{\mathcal{L}})$ be an expanded MDP and $A_{\varphi}$$=$$(\mathcal{Q}, q_0, 2^{\mathcal{AP}}, \delta_{\varphi}, \mathcal{F})$ be a DFA. The \textit{product MDP} $\mathcal{M}_p$$=$$(S_p, s_{0_p}, \mathcal{A}, \mathbb{P},\mathcal{AP}, \mathcal{L}_p, \mathcal{F}_p)$ is a tuple where
\begin{itemize}
\item $S_p$$=$$\overline{S} $$\times$$ Q$,
\item $s_{0_p}=(\overline{s}_0,q)$ such that $q=\delta(q_0,\overline{\mathcal{L}}(\overline{s}_0))$,
\item $\mathbb{P}((s,q), a, (s',q'))$=$\begin{cases} \overline{\mathcal{P}}_{s,a,s'} & \text{if} \quad q'=\delta(q,\overline{\mathcal{L}}(s')) \\ 0 & \text{otherwise}, \end{cases}$
\item $\mathcal{L}_p((s,q))=\{q\}$,
\item $\mathcal{F}_p$$=$$\overline{S}\times \mathcal{F}$.
\end{itemize}
\end{definition}}

The incentive design problem \eqref{main_prob_1}-\eqref{maxreach} can now be solved on the product MDP $\mathcal{M}_p$ in three steps. First, we partition the states of $\mathcal{M}_p$ into three disjoint sets. Let $B$$:=$$\mathcal{F}_p$, $S_0$ be the set of states that have zero probability of reaching the set $B$, and $S_r$$:=$$S_p\backslash B\cup S_0$. Second, we form the modified product MDP $\mathcal{M}_p'$ by making all states $B\cup S_0$ absorbing. Finally, we apply the methods developed in Section \ref{thecostofcontrol} to synthesize an $\overline{\epsilon}$-optimal incentive sequence on $\mathcal{M}_p'$.

 Note that the incentive sequence is designed on the product MDP $\mathcal{M}_p$. Therefore, the principal must share the DFA structure, i.e., it's objective, with the agent to be able to use the computed design. However, for the existence of a solution to the problem \eqref{main_prob_1}-\eqref{maxreach}, the incentive sequence should be designed on the MDP $\mathcal{M}$. The following example illustrates that the problem \eqref{main_prob_1}-\eqref{maxreach} may have no feasible solution, even though the existence of an $\overline{\epsilon}$-optimal incentive sequence on $\mathcal{M}_p$ is guaranteed.
 
 \begin{figure}[b!]\centering\vspace{-0.4cm}
\scalebox{0.8}{\begin{tikzpicture}[->, >=stealth', auto, semithick, node distance=2cm]

    \tikzstyle{every state}=[fill=white,draw=black,thick,text=black,scale=0.7]

    \node[state,initial below,initial text=] (s_0) {$s_0,A$};
    \node[state] (s_1) [below right =3mm and 15mm  of s_0]  {$s_1, A$};
    \node[state] (s_2) [above right =2mm and 25mm of s_0]  {$s_2, B$};
    \node[state] (s_3) [left=15mm of s_0]  {$s_3, C$};

\path
(s_0)  edge[out=70, in=110, loop]    node[label={[xshift=0cm, yshift=0cm]$a_1,0.2$}]{}     (s_0)
(s_0)	 edge     node[label={[xshift=-0.3cm, yshift=-0.9cm]$a_1, 0.4$}]{}     (s_1)
(s_0)	 edge[in=220, out=0]   node[label={[xshift=0.3cm, yshift=-0.7cm]$a_1, 0.4$}]{}     (s_2)
(s_0)	 edge   node[label={[xshift=0cm, yshift=0cm]$a_2, 1$}]{}     (s_3)
(s_1)  edge[out=-20, in=20, loop]    node[label={[xshift=0.5cm, yshift=-0.4cm]$a_1,1$}]{}     (s_1)
(s_3)  edge[out=160, in=200, loop]    node[label={[xshift=-0.5cm, yshift=-0.4cm]$a_1,1$}]{}     (s_3)

(s_2)  edge[out=170, in=40]    node[label={[xshift=0cm, yshift=0.1cm]$a_1,1$}]{}     (s_0);
\end{tikzpicture}}
\caption{An MDP example for which there exists no feasible incentive design for the scLTL specification $\varphi$$=$$\lozenge(B$$\land$$\lozenge C)$.}\label{impossibility}
\end{figure}
 \noindent\textbf{Example 2:} Consider the MDP given in Fig. \ref{impossibility}, where the numbers next to actions $a_i$ represent the transition probabilities, e.g., $\mathcal{P}_{s_0,a_1,s_1}$$=$$0.4$, and the letters next to state numbers represent labels, e.g., $\mathcal{L}(s_0)$$=$$A$. Let the agent's decision horizon be $N$$=$$3$, and the reward function $\mathcal{R}$ be such that $\mathcal{R}(s_0,a_1)$$=$$1$ and $\mathcal{R}(s,a)$$=$$0$ otherwise.
 Additionally, let the principal's objective be expressed by the scLTL formula $\varphi$$=$$\lozenge(B$$\land$$\lozenge C)$, i.e., first visit state $B$ and then state $C$. The maximum probability of satisfying $\varphi$ is $x^{\star}_0$$=$$0.5$, which can be computed by solving an LP \cite{Model_checking}. The value $x^{\star}_0$ is attainable if and only if the agent takes the action $a_2$$\in$$\mathcal{A}(s_0)$ with probability 1 after visiting state $s_2$. 
 
 The principal should decide on which actions to incentivize in the first three stages $t$$=$$1,2,3$ since the agent's decision horizon is $N$$=$$3$. Clearly, the action $a_1$ should be incentivized for $t$$=$$1,2$ so that the agent visits state $s_2$. At $t$$=$$3$, the agent will be in state $s_0$ with nonzero probability. Now, if the principal incentivize $a_1$, the agent will take action $a_2$$\in$$\mathcal{A}(s_0)$ with probability less then 1 after visiting state $s_2$. On the other hand, if $a_2$ is incentivized, then the agent cannot satisfy the specification with probability higher than 0.46. Consequently, no incentive design on the given MDP can guarantee the satisfaction of the specification $\varphi$ with maximum probability. $\triangleleft$

 We now present a sufficient condition on the structure of the MDP $\mathcal{M}$ which guarantees the existence of an $\overline{\epsilon}$-optimal incentive design on $\mathcal{M}$.

 For the product MDP $\mathcal{M}_p$ and a policy $\pi$$\in$$\Pi(\mathcal{M}_p)$, let $M^{\pi}_{s,t}$$:=$$\{q$$\in$$\mathcal{Q}$$:$$\sum_{a\in\mathcal{A}}\mu^{\pi}_t((s,q),a)$$>$$0 \}$ be the set of occupied memory states when the agent is in state $s$$\in$$\overline{S}$ at stage $t$$\in$$\mathbb{N}$.
{\setlength{\parindent}{0cm}
 \noindent \begin{thm} \label{theorem_suff} For an MDP $\mathcal{M}$ and a decision horizon $N$, let $\overline{S}$ be the finite set of states for the expanded MDP. There exists an $\overline{\epsilon}$-optimal incentive design on $\mathcal{M}$ if there exists an $\overline{\epsilon}$-optimal incentive design on $\mathcal{M}_p$ such that the agent's optimal policy $\pi$$\in$$\Pi(\mathcal{M}_p)$ under the provided incentives satisfies $\lvert M^{\pi}_{s,t}\rvert$$\leq$$1$ for all $s$$\in$$\overline{S}$ and $t$$\in$$\mathbb{N}$.
 \end{thm}}
 \noindent\textbf{Proof:}
 Let $\Gamma$ be an $\overline{\epsilon}$-optimal incentive design on $\mathcal{M}_p$ with the desired property. Note that the function in \eqref{optincentives} is a mapping from the expanded MDP $\overline{\mathcal{M}}$ to the MDP $\mathcal{M}$ that preserves $\overline{\epsilon}$-optimality of the incentive design. Therefore, in what follows, we construct an incentive mapping from $\mathcal{M}_p$ to $\overline{\mathcal{M}}$ that preserves $\overline{\epsilon}$-optimality of the incentive design $\Gamma$, and conclude the result.
 
An $\overline{\epsilon}$-optimality-preserving mapping $\psi$ such that $\Gamma'$$=$$\{\gamma'_1,\gamma'_2,\ldots\}$$:=$$\psi(\Gamma)$ where $\Gamma$$=$$\{\gamma_1,\gamma_2,\ldots\}$ is given as follows. For a given $t$$\in$$\mathbb{N}$,
\begin{itemize}
    \item if $\lvert M^{\pi}_{s,t}\rvert$$=$$0$, $\gamma'_t(s,a)$$:=$$\gamma_t((s,q),a)$ for an arbitrary $q$$\in$$\mathcal{Q}$ and for all $a$$\in$$\mathcal{A}(s)$,
    \item if $\lvert M^{\pi}_{s,t}\rvert$$=$$1$, $\gamma'_t(s,a)$$:=$$\gamma_t((s,q),a)$ for $q$$\in$$M^{\pi}_{s,t}$ and for all $a$$\in$$\mathcal{A}(s)$. $\Box$
\end{itemize}

The following corollary follows from the fact that the principal can induce a stationary deterministic agent policy on the product MDP through the methods explained in Section \ref{computational_section}.
{\setlength{\parindent}{0cm}
 \noindent \begin{Corollary} For an MDP $\mathcal{M}$, there exists an $\overline{\epsilon}$-optimal incentive design if  $\mathcal{P}_{s,a,s'}$$\in$$\{0,1\}$ for all $s,s'$$\in$$S$ and $a$$\in$$\mathcal{A}(s)$.
 \end{Corollary}}

 Finally, we provide a sufficient condition on the agent's decision horizon $N$ that ensures the existence of an $\overline{\epsilon}$-optimal incentive design on $\mathcal{M}$.
 {\setlength{\parindent}{0cm}
 \noindent \begin{prop} \label{N1prop}For an MDP $\mathcal{M}$, there exists an $\overline{\epsilon}$-optimal incentive design if the agent's decision horizon is $N$$=$$1$.
 \end{prop}}
 \noindent\textbf{Proof (Sketch):} There is a one-to-one correspondence between the paths of the product MDP $\mathcal{M}_p$ and the MDP $\mathcal{M}$ \cite{Model_checking}. Therefore, the principal can observe the path followed by the agent on $\mathcal{M}$, and provide the incentives according to the corresponding path on $\mathcal{M}_p$ at each stage. Because $N$$=$$1$, the principal knows the memory state occupied by the agent at each stage. Consequently, it becomes possible to map the incentives from $\mathcal{M}_p$ to $\mathcal{M}$ at each stage. $\Box$
 
 \section{Numerical Simulations}
 In this section, we demonstrate the proposed incentive design methods on two simple motion planning examples. Considering the availability of off-the-shelf solvers, e.g., Gurobi\cite{gurobi}, MOSEK\cite{mosek}, that can efficiently solve large-scale linear optimization problems, we restrict our attention to small scale examples to better emphasize the properties of the proposed methods. We synthesize the incentive sequences for the following examples through the use of MOSEK \cite{mosek} solver together with CVXPY \cite{cvxpy} interface. 
 \subsection{Incentives for reachability objectives}
 In this example, we consider a $5$$\times$$5$ grid world environment, shown in Fig. \ref{first_example}, and an agent with decision horizon $N$$=$$1$. At each state, the agent has four actions, i.e., $\mathcal{A}$$=$$\{left,right,up,down\}$, and a transition to the chosen direction occurs with probability 1. If the adjacent state in the chosen direction is the boundary of the environment, the agent stays in its current state. A reward function for the agent is generated by choosing all rewards $\mathcal{R}(s,a)$ from the set $\{0,1,\ldots,9\}$ uniformly randomly.

 The agent starts from the bottom left corner, i.e., \textit{Start} state in Fig. \ref{first_example}, and aims to maximize its immediate reward at each stage. The principal provides incentives to the agent so that the agent reaches the top right corner, i.e., \textit{Target} state in Fig. \ref{first_example}. 

\begin{figure}[b]\vspace{-0.3cm}\centering 
\scalebox{0.2}{
\begin{tikzpicture}
\draw[black,line width=1pt] (0,0) grid[step=4] (20,20);
\draw[black,line width=4pt] (0,0) rectangle (20,20);

\node at (1.5,1)  { \textbf{\fontsize{30pt}{30pt}\selectfont Start}};
\node at (18,19)  { \textbf{\fontsize{30pt}{30pt}\selectfont Target}};

\node (1) at (2,2) {};
\node (2) at (6,2) {};

\path[color=blue, ultra thick,decoration={
    markings,
    mark=at position 0.6 with {\arrow[scale=4,>=stealth]{>}}}]
(1) edge [postaction={decorate}, bend right=60pt] node {} (2)
(2) edge [postaction={decorate},bend right=60pt] node {} (1);

\draw [color=red,ultra thick, decoration={markings,mark=at position 0.4 with {\arrow[scale=4,>=stealth]{>}}},postaction={decorate}] (2,2) -- (6,2);

\draw [color=red,ultra thick, decoration={markings,mark=at position 0.4 with {\arrow[scale=4,>=stealth]{>}}},postaction={decorate}] (6,2) -- (10,2);

\draw [color=red,ultra thick, decoration={markings,mark=at position 0.4 with {\arrow[scale=4,>=stealth]{>}}},postaction={decorate}] (10,2) -- (10,6);

\draw [color=red,ultra thick, decoration={markings,mark=at position 0.4 with {\arrow[scale=4,>=stealth]{>}}},postaction={decorate}] (10,6) -- (14,6);

\draw [color=red,ultra thick, decoration={markings,mark=at position 0.4 with {\arrow[scale=4,>=stealth]{>}}},postaction={decorate}] (14,6) -- (14,10);

\draw [color=red,ultra thick, decoration={markings,mark=at position 0.4 with {\arrow[scale=4,>=stealth]{>}}},postaction={decorate}] (14,10) -- (18,10);

\draw [color=red,ultra thick, decoration={markings,mark=at position 0.4 with {\arrow[scale=4,>=stealth]{>}}},postaction={decorate}] (18,10) -- (18,14);

\draw [color=red,ultra thick, decoration={markings,mark=at position 0.4 with {\arrow[scale=4,>=stealth]{>}}},postaction={decorate}] (18,14) -- (14,14);

\draw [color=red,ultra thick, decoration={markings,mark=at position 0.4 with {\arrow[scale=4,>=stealth]{>}}},postaction={decorate}] (14,14) -- (14,18);

\draw [color=red,ultra thick, decoration={markings,mark=at position 0.4 with {\arrow[scale=4,>=stealth]{>}}},postaction={decorate}] (14,18) -- (18,18);

\end{tikzpicture}
}

\caption{The motion of an agent on a grid world. The agent's decision horizon is $N$$=$$1$, and it starts from the \textit{Start} state. The principal's objective is to induce an agent policy that reaches the \textit{Target} state with probability 1. Blue arrows indicate the agent's optimal policy in the absence of incentives, and red arrows indicate the agent's optimal policy under the provided incentives. }
\label{first_example}
\end{figure}

 In the absence of incentives, i.e., $\gamma_t(s,a)$$=$$0$ for all $t$$\in$$\mathbb{N}$, the agent's optimal path is shown by blue arrows in Fig. \ref{first_example}. Under its optimal policy, the agent cycles between two states infinitely often. Through the methods explained in Section \ref{thecostofcontrol}-\ref{computational_section}, we synthesize an incentive sequence for the agent so that it reaches the target state with probability 1. The agent's optimal path under the provided incentives is shown by red arrows in Fig. \ref{first_example}. The total cost to the principal is computed as $9$$+$$10\epsilon$ units of resources (UR) where $\epsilon$$>$$0$ is an arbitrarily small constant.
 
 As can be seen from Fig. \ref{first_example}, under the provided incentives, the agent follows the lowest cost path rather than the shortest one to the target state. Specifically, the shortest path would take $8$ stages to reach the target state and cost $12$ UR to the principal, whereas the lowest cost path takes $10$ stages to reach the target and cost $9$ UR. Quantitatively, the proposed incentive design allows the principal to save $25\%$ of the resources that would be paid to the agent if it was to follow the shortest path.
\subsection{Incentives for general scLTL specifications}
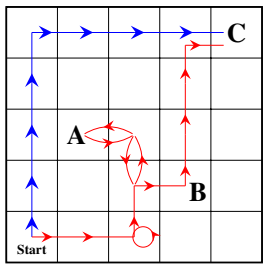
\begin{figure}[H]\vspace{-0.2cm}\centering
\scalebox{0.17}{
\begin{tikzpicture}
\draw[black,line width=1pt] (0,0) grid[step=4] (20,20);
\draw[black,line width=4pt] (0,0) rectangle (20,20);

\node at (2,1)  { \textbf{\fontsize{30pt}{30pt}\selectfont Start}};
\node at (5.5,10)  { \textbf{\fontsize{60pt}{60pt}\selectfont A}};
\node at (15,5.5) { \textbf{\fontsize{60pt}{60pt}\selectfont B}};
\node at (18,18) { \textbf{\fontsize{60pt}{60pt}\selectfont C}};
\node (1) at (10,10) {};
\node (2) at (6,10) {};
\node (3) at (10,6) {};
\draw [color=blue,ultra thick, decoration={markings,mark=at position 0.4 with {\arrow[scale=5,>=stealth]{>}}},postaction={decorate}] (2,2) -- (2,6);
\draw [color=blue,ultra thick, decoration={markings,mark=at position 0.2 with {\arrow[scale=5,>=stealth]{>}}},postaction={decorate}] (2,6) -- (2,10);
\draw [color=blue,ultra thick, decoration={markings,mark=at position 0.2 with {\arrow[scale=5,>=stealth]{>}}},postaction={decorate}] (2,10) -- (2,14);
\draw [color=blue,ultra thick, decoration={markings,mark=at position 0.2 with {\arrow[scale=5,>=stealth]{>}}},postaction={decorate}] (2,14) -- (2,18);

\draw [color=blue,ultra thick, decoration={markings,mark=at position 0.4 with {\arrow[scale=5,>=stealth]{>}}},postaction={decorate}] (2,18) -- (6,18);
\draw [color=blue,ultra thick, decoration={markings,mark=at position 0.2 with {\arrow[scale=5,>=stealth]{>}}},postaction={decorate}] (6,18) -- (10,18);
\draw [color=blue,ultra thick, decoration={markings,mark=at position 0.2 with {\arrow[scale=5,>=stealth]{>}}},postaction={decorate}] (10,18) -- (14,18);
\draw [color=blue,ultra thick, decoration={markings,mark=at position 0.2 with {\arrow[scale=5,>=stealth]{>}}},postaction={decorate}] (14,18) -- (17,18);

\draw [color=red,ultra thick, decoration={markings,mark=at position 0.4 with {\arrow[scale=4,>=stealth]{>}}},postaction={decorate}] (2,2) -- (6,2);

\draw [color=red,ultra thick, decoration={markings,mark=at position 0.2 with {\arrow[scale=4,>=stealth]{>}}},postaction={decorate}] (6,2) -- (10,2);

\draw [color=red,ultra thick, decoration={markings,mark=at position 0.1 with {\arrow[scale=3.5,>=stealth]{>}}},postaction={decorate}] (10.7,2) circle (0.8cm);
\draw [color=red,ultra thick, decoration={markings,mark=at position 0.4 with {\arrow[scale=4,>=stealth]{>}}},postaction={decorate}] (10,2) -- (10,6);

\draw [color=red,ultra thick, decoration={markings,mark=at position 0.4 with {\arrow[scale=4,>=stealth]{>}}},postaction={decorate}] (10,6) -- (14,6);

\draw [color=red,ultra thick, decoration={markings,mark=at position 0.4 with {\arrow[scale=4,>=stealth]{>}}},postaction={decorate}] (14,6) -- (14,10);
\draw [color=red,ultra thick, decoration={markings,mark=at position 0.4 with {\arrow[scale=4,>=stealth]{>}}},postaction={decorate}] (14,10) -- (14,14);

\draw [color=red,ultra thick, decoration={markings,mark=at position 0.4 with {\arrow[scale=4,>=stealth]{>}}},postaction={decorate}] (14,14) -- (14,17);

\draw [color=red,ultra thick, decoration={markings,mark=at position 0.4 with {\arrow[scale=4,>=stealth]{>}}},postaction={decorate}] (14,17) -- (17,17);

\path[color=red, ultra thick,decoration={
    markings,
    mark=at position 0.6 with {\arrow[scale=4,>=stealth]{>}}}]
(1) edge [postaction={decorate}, bend right] node {} (2)
(2) edge [postaction={decorate},bend right] node {} (1);

\path[color=red, ultra thick,decoration={
    markings,
    mark=at position 0.6 with {\arrow[scale=4,>=stealth]{>}}}]
(1) edge [postaction={decorate}, bend right] node {} (3)
(3) edge [postaction={decorate},bend right] node {} (1);

\end{tikzpicture}
}
\caption{The motion of an agent on a grid world. The agent's decision horizon is $N$$=$$4$, and it starts from the \textit{Start} state. The principal's objective is to induce an agent policy that satisfies the scLTL specification $\varphi$$=$$\lozenge(A \land \lozenge( B \land \lozenge C))$, i.e., first visit $A$, then $B$, and then $C$. The optimal path of the agent in the absence of incentives is shown by blue arrows (top path). Red arrows indicate the agent's optimal path under the provided incentives (bottom path). }
\label{grid_graph}
\end{figure}

In this example, we consider the same grid world environment introduced in the previous example with different state labels. The agent's decision horizon is $N$$=$$4$, and its objective is to reach the state labeled as $C$ in Fig. \ref{grid_graph}. The principal's objective is to induce an agent policy that satisfies the scLTL specification $\varphi$$=$$\lozenge(A \land \lozenge( B \land \lozenge C))$, i.e., the agent should first visit state $A$, then $B$, and then $C$, with probability 1. 

The agent receives the reward of $2$ for transitioning to the top left state and the reward of $5$ for transitioning to the top right state. Its optimal path in the absence of incentives is shown in Fig. \ref{grid_graph} with blue arrows (top path). We synthesize an optimal incentive sequence under which the agent's optimal path is shown in Fig. \ref{grid_graph} with red arrows (bottom path).

The total cost of the incentives to the principal is computed as $2$$+$$13\epsilon$ units of resources. Specifically, the principal provides $2$$+$$\epsilon$ incentives for the $right$ action in the start state and then $\epsilon$ incentives at each stage for desired actions. An interesting property of the incentivized (red) path is that the agent stays in the same state in third stage by taking $down$ action. This is due to the fact that the state $s$ on the left of the state labeled as $B$ has value $V_n(s)$$=$$0$ for all $n$. Therefore, the principal wants that state to be the agent's initial state when it computes its second $4$-stage policy. By doing so, the principal ensures that the states $s'$ occupied by the agent in the next 4 stages will always have a value zero, i.e., $V_n(s')$$=$$0$ if $\sum_{a\in\mathcal{A}(s')}\mu_{4+n}^{\pi}(s',a)$$>$$0$, and therefore the cost of control will only be $\epsilon$.

\section{Conclusions and Future Directions}

We considered a principal-agent model and studied the problem of designing an optimal sequence of incentives that the principal should offer to the agent in order to induce a desired agent behavior expressed as a syntactically co-safe linear temporal logic (scLTL) formula. For reachability objectives, we presented a polynomial-time algorithm to synthesize an incentive design that minimizes the cost to the principal. By providing an example scenario, we showed that a feasible incentive design may not exists for general scLTL formulae, and the principal may need to share its objective with the agent to induce the desired behavior. Furthermore, we provided sufficient conditions under which the principal can induce the desired behavior without sharing the scLTL formula with the agent. 

The results that we present in this paper are obtained under the assumptions that the agent's reward function and the length of its decision horizon are known by the principal. An interesting future direction may be to develop methods to infer the length of the agent's decision horizon through perfect/imperfect observations, or to design an incentive sequence that does not require the knowledge of the length of the decision horizon.

  \bibliographystyle{IEEEtran}
\bibliography{main.bib}
\newpage
\appendices \section{Proof of Proposition \ref{comp_prop}}\label{proof_appendix}

Note that any deterministic policy constructed from the optimal decision variables $\lambda^{\star}(s,a)$ of  \eqref{obj22}-\eqref{const_last2} can only violate the reachability constraint \eqref{cons_1}. In other words, the constructed policy is guaranteed to minimize the expected total cost.

We will need the following result to prove Proposition \ref{comp_prop}. For a given policy $\pi$, let $Reach^{\pi}(s,s')$ denote the probability of reaching $s'$ from $s$ under $\pi$. Note that $Reach^{\pi}(s,s)$$=$$\sum_{a\in\mathcal{A}(s)}d(s,a)Reach^{\pi}((s,a),s)$ where $Reach^{\pi}((s,a),s)$ is the probability of reaching state $s$ from state action pair $(s,a)$ under the policy $\pi$. Finally, let $\xi^{\pi}(s)$$:=$$\sum_{a\in\mathcal{A}(s)}\xi^{\pi}(s,a)$, and note that for any $\xi^{\pi}(s)$$<$$\infty$, we have \cite{serfozo}
\begin{align}
    \xi^{\pi}(s)=\frac{Reach^{\pi}(s_0,s)}{1-\sum_{a\in\mathcal{A}(s)}d(s)(a)Reach^{\pi}((s,a),s)}.
\end{align}

To prove the claim of Proposition \ref{comp_prop}, we show that any policy that is constructed by choosing actions $a$$\in$$\mathcal{A}(s)$ such that $\lambda^{\star}(s,a)$$>$$0$ deterministically is optimal. 

Let $\overline{\pi}$$=$$\{\overline{d},\overline{d},\ldots\}$ be the stationary randomized policy constructed from $\lambda^{\star}(s,a)$ of LP \eqref{opt_2} through the formula \eqref{classic_randomized}. Additionally, let $\widetilde{\pi}$$=$$\{\widetilde{d},\widetilde{d},\ldots\}$ be a stationary randomized policy such that $\widetilde{d}(s^{\star},a^{\star})$$=$$1$, and $\widetilde{d}(s)$$=$$\overline{d}(s)$ for all $s$$\in$$S\backslash \{s^{\star}\}$. Informally, in state $s^{\star}$, we choose one of the \textit{active} actions deterministically and do not change the rest of the policy.

We first show that $\text{Pr}^{\widetilde{\pi}}(s_0$$\models$$ \varphi)$$<$$x^{\star}_{s_0}$ implies $Reach^{\widetilde{\pi}}((s^{\star},a^{\star}),s^{\star})$$=$$1$. Then by showing that $Reach^{\widetilde{\pi}}((s^{\star},a^{\star}),s^{\star})$$=$$1$ cannot be true, we conclude that $\text{Pr}^{\widetilde{\pi}}(s_0$$\models$$ \varphi)$$=$$x^{\star}_{s_0}$.

As for the first claim, suppose for contradiction that $\text{Pr}^{\widetilde{\pi}}(s_0$$\models$$ \varphi)$$<$$x^{\star}_{s_0}$ and $Reach^{\widetilde{\pi}}((s^{\star},a^{\star}),s^{\star})$$<$$1$. Note that if 
$Reach^{\widetilde{\pi}}((s^{\star},a^{\star}),s^{\star})$$<$$1$, then $Reach^{\widetilde{\pi}}(s^{\star},s^{\star})$$<$$1$. Therefore, $Reach^{\widetilde{\pi}}(s,s)$$<$$1$ for all $s$$\in$$S_r$ satisfying $Reach^{\widetilde{\pi}}(s^{\star},s)$$>$$0$. Additionally, as $\overline{d}(s')$$=$$\widetilde{d}(s')$ for all $s'$ such that $Reach^{\widetilde{\pi}}(s^{\star},s')$$=$$0$, we have $Reach^{\widetilde{\pi}}(s',s')$$<$$1$. Consequently, probability of leaving the set $S_r$ is 1. Since all actions that are chosen by policy $\widetilde{\pi}$ satisfy $x^{\star}_s$$=$$\mathcal{P}_{s,a,s'}x^{\star}_{s'}$
where $x^{\star}_s$ is the maximum probability of reaching the set $B$ from the state $s$ (see e.g. Chapter 10 in \cite{Model_checking}), 
probability of entering the set $B$ must be equal to $x^{\star}_{s_0}$. This raises a contradiction.

As for the second claim, suppose that $Reach^{\widetilde{\pi}}((s^{\star},a^{\star}),s^{\star})$$=$$1$. Then, $Reach^{\overline{\pi}}((s^{\star},a^{\star}),s^{\star})$$=$$1$ since $\overline{\pi}$ differs from $\widetilde{\pi}$ only in the state $s^{\star}$. We now construct a policy $\hat{\pi}$ such that $\hat{d}(s)$$=$$\overline{d}(s)$ for all $s$$\in$$S\backslash \{s^{\star}\}$, $\hat{d}(s^{\star},a^{\star})$$=$$0$, and
\begin{align}\label{cool}
    \hat{d}(s^{\star},a_i)=\frac{\overline{d}(s^{\star},a_i)}{\sum_{a\in\mathcal{A}(s^{\star})\backslash\{a^{\star}\}}\overline{d}(s^{\star},a_i)}.
 \end{align}
Note that $\hat{\pi}$ satisfies $\text{Pr}^{\hat{\pi}}(s_0$$\models$$ \varphi)$$=$$x^{\star}_{s_0}$. By showing that $\hat{\pi}$ attains an objective value in \eqref{opt_2} that is strictly smaller than the policy $\overline{\pi}$, we will conclude that $Reach^{\widetilde{\pi}}((s^{\star},a^{\star}),s^{\star})$$=$$1$ cannot be possible.

For the ease of notation, let $\overline{a}_i$$:=$$\overline{d}(s^{\star},a_i)$, $\overline{R}_i$$:=$$Reach^{\overline{\pi}}((s^{\star},a_i),s^{\star})$, and $\hat{a}_i$$:=$$\hat{d}(s^{\star},a_i)$, and $\hat{R}_i$$:=$$Reach^{\hat{\pi}}((s^{\star},a_i),s^{\star})$. Without loss of generality, we choose $a_1$$=$$a^{\star}$.
By the construction of $\hat{\pi}$, it can be shown that
\begin{align}
    \xi^{\hat{\pi}}(s^{\star})=(1-\overline{a}_1)\xi^{\overline{\pi}}(s^{\star})-\frac{\overline{a}_1(\overline{R}_1-1)(1-\overline{a}_1)}{C}
\end{align}
where $C$$:=$$(1-\sum_{i=1}^n \overline{a}_i\overline{R}_i)(1-a_1-\sum_{i=2}^n \overline{a}_i\overline{R}_i)$$>$$0$. Note that $\xi^{\hat{\pi}}(s^{\star},a_i)(1-\overline{a}_1)$$=$$\xi^{\overline{\pi}}(s^{\star},a_i)$ due to \eqref{cool}. Then, since $\overline{a}_1$$>$$0$ and $\overline{R}_1$$=$$1$, we have
\begin{align}
    \xi^{\hat{\pi}}(s^{\star},a_i)\leq \xi^{\overline{\pi}}(s^{\star},a_i)
\end{align}
for all $a_i$ $i$$=$$2,3,\ldots,n$ and $\xi^{\hat{\pi}}(s^{\star},a_1)$$<$$\xi^{\overline{\pi}}(s^{\star},a_1)$. Consequently, $\hat{\pi}$ attains an objective value in \eqref{opt_2} that is strictly smaller than the policy $\overline{\pi}$. 

Finally, since $Reach^{\overline{\pi}}((s^{\star},a^{\star}),s^{\star})$$=$$1$ cannot be true, $Reach^{\widetilde{\pi}}((s^{\star},a^{\star}),s^{\star})$$=$$1$ cannot be true. If $Reach^{\widetilde{\pi}}((s^{\star},a^{\star}),s^{\star})$$=$$1$ is not true, $\text{Pr}^{\widetilde{\pi}(s_0\models\varphi)}$$<$$x^{\star}_{s_0}$ is not true. This concludes the proof. $\Box$
\end{document}